\documentclass[a4article,12pt]{article}

\usepackage[T2A]{fontenc}
\usepackage[russian,english]{babel}
\usepackage[cp1251]{inputenc}
\usepackage{amsthm,amsmath}

\usepackage{amssymb}
\usepackage{graphicx}
\usepackage{color}

\makeatletter
\renewcommand{\@biblabel}[1]{#1.}
\makeatother

\theoremstyle{plain}
\newtheorem{theorem}{Теорема}
\newtheorem{lemma}{Лемма}

\theoremstyle{definition}

\newtheorem{corollary}{Следствие}

\begin{document}
\selectlanguage{russian}
\begin{center}
\textbf{Об индексе дефекта векторного оператора Штурма--Лиувилля}\footnote{Работа будет опубликована в феврале 2016 г. в журнале "Математические заметки" (том 99, вып.2).}\\
К.А.Мирзоев, Т.А. Сафонова\footnote{Работа первого автора выполнена при финансовой поддержке РНФ (грант № 14-11-00754).
Второй автор поддержан Минобрнауки РФ (грант Президента РФ № МК-3941.2015.1), РФФИ (гранты №№ 14-01-31136-мол а, 14-01-00349, 15-31-50259 ) и фондом содействия отечественной науки.}
\end{center}

{\bf 1. Введение. Предварительные сведения}

{ 1.1.} Пусть $R_+:=[0,+\infty)$ и пусть матриц-функции $P, Q$ и $R$ порядка $n$ ($n\in\mathcal N$), определённые на полуоси $R_+$, таковы, что $P(x)$ - невырожденная, $P(x)$ и $Q(x)$ - эрмитовы матрицы при $x\in R_+$, а элементы матриц-функций $P^{-1}$, $Q$ и $R$ измеримы на $R_+$ и суммируемы на каждом её замкнутом конечном подынтервале. Определим первую квазипроизводную заданной локально абсолютно непрерывной вектор-функции $f=(f_1(x),f_2(x),$ $\ldots,f_n(x))^t$ (здесь и далее везде, $t$ - символ транспонирования), полагая
$f^{[1]}:=P(f^{\prime}-Rf).$ Предположим, что функция $f^{[1]}$ также локально абсолютно непрерывна, и определим вторую квазипроизводную $f^{[2]}$, полагая
$f^{[2]}:=(f^{[1]}) ^{\prime}+R^*f^{[1]}-Qf,$
где $^*$ - символ сопряжения, и квазидифференциальное выражение, полагая
$l[f](x):=-f^{[2]}(x).$ Таким образом,
\begin{equation}
\label{trivial}
l[f]=-(P(f^{\prime}-Rf))^{\prime}-R^*P(f^{\prime}-Rf)+Qf,
\end{equation}
а область определения $\Delta$ выражения  $l[f]$ --- это множество всех комплекснозначных вектор-функций $f$, таких что $f$ и $f^{[1]}$ локально абсолютно непрерывны на $R_+$ и при
$f\in \Delta$ выражение $l[f]$ $(\in \mathcal{L}^1_{loc}(R_+))$ определяется п.в. по формуле~(\ref{trivial}).
Кроме того, для любых двух вектор-функций $f,g \in \Delta$ справедлива следующая лемма - векторный аналог тождества Грина:
\begin{lemma}
\label {formulaGrina}
Пусть $P(x)$, $Q(x)$ и $R(x)$ - матриц-функции порядка $n$, удовлетворяющие перечисленным выше условиям. Тогда для любых двух вектор-функций $u,v\in \Delta$ и для любых двух чисел $\alpha$ и $\beta$ таких, что $0 \le \alpha \le \beta < \infty$, справедлива формула
\begin{equation}
 \label{1.1}
  \int\limits_{\alpha}^{\beta} \{(l[u](x),v(x))  - (u(x),l[v](x))\}dx  = [u(x),v(x)] (\beta)-[u(x),v(x)](\alpha),
 \end{equation}

где  $(g,h)=\sum\limits_{s=1}^{n}g_s\overline{h_s}$ - скалярное произведение векторов $g$ и $h$, а форма $[u,v]$ определяется равенством
 $ [u,v] (x):=(u^{[1]}(x),v(x))-(u(x),v^{[1]}(x))$.
 \end{lemma}

Пусть далее $\mathcal{L}^2_n(R_+)$ -- пространство классов эквивалентности всех комплекснозначных измеримых вектор-функций $f$, у которых сумма квадратов модулей компонент интегрируема по Лебегу на $R_+$. В литературе, посвященной спектральной теории обыкновенных дифференциальных операторов, хорошо известна процедура, с помощью которой определяются минимальный и максимальный операторы
$L_0$ и $L_1$ соответственно, порожденные выражением  $l[f]$ в гильбертовом пространстве $\mathcal{L}^2_n(R_+)$.
Операторы $L_0$, $L_1$ и операторы, связанные с ними, называются векторными операторами Штурма -- Лиувилля. Используя формулу (\ref{1.1}), можно доказать, что оператор $L_0$ является замкнутым симметрическим оператором. Кроме того, область определения оператора $L_0$ всюду плотна в $\mathcal{L}^2_n(R_+)$. Пусть пара $(d_+,d_-)$ - индекс дефекта оператора $L_0$.  Согласно \cite{Naim}-\cite{Anderson},  дефектные числа $d_+$ и $d_-$ оператора $L_0$ совпадают с максимальным числом линейно независимых решений уравнения
\begin{equation}
\label{000}
l[f]=\lambda f,
 \end{equation}
 принадлежащих пространству $\mathcal{L}^2_n(R_+)$, когда параметр $\lambda$ берётся из верхней $(\Im\lambda > 0)$ или нижней $(\Im \lambda < 0)$ полуплоскости соответственно, удовлетворяют неравенствам $n \le d_+,d_- \le 2n$ и, кроме того,  $d_+=2n$ тогда и только тогда, когда $d_-=2n$. Случай $d_+=d_-=2n$ реализуется тогда и только тогда, когда все решения уравнения (\ref{000}) при всех $\lambda\in C$ принадлежат пространству $\mathcal{L}^2_n(R_+)$.
Используя аналогию со спектральной теорией скалярных операторов Штурма-Лиувилля на полуоси, иногда говорят, что для выражения $l[f]$ (оператора $L_0$) имеет место случай предельной точки, если $d_+=d_-=n$, если же $d_+=d_-=2n$, то говорят, что для выражения $l[f]$ имеет место случай предельного круга (см., напр., \cite{Anderson}).

Уравнение (\ref{000}) равносильно системе дифференциальных уравнений первого порядка
\begin{equation}
\label{ss}
Y^{\prime}=(F-\Lambda)Y,
\end{equation}
где $Y=(f,f^{[1]})^t$, матрицы $F$ и $\Lambda$ порядка $2n$ имеют вид
$$F=\begin{pmatrix}
R & P^{-1}\\
Q& {-R^*}\\
\end{pmatrix},\,\,\,\,\,\,
\Lambda=\begin{pmatrix}
O & O\\
\lambda I & O\\
\end{pmatrix},$$
а $O$ и $I$, как обычно, - нулевая и единичная матрицы порядка $n$ соответственно.
Равносильность этих уравнений понимается в том смысле,
  что если $n$-компонентная вектор-функция $f$ является решением  (\ref {000}), то $2n$-компонентная вектор-функция $Y=(f,f^{[1]})^t$ является решением системы (\ref {ss}) и наоборот, если   $Y$ - решение этой системы, то вектор-функция $f$, составленная из его первых $n$-компонент, - решение системы уравнений   ~(\ref {000}).

Заметим, что условия на элементы матриц $P,\,Q$ и $R$, перечисленные выше,  обеспечивают для системы (\ref {ss}) справедливость теоремы существования и единственности решения задачи Коши, поставленной в произвольной точке $R_+$.

Используя терминологию из теории линейных квазидифференциальных уравнений, говорят, что квазипроизводные $f^{[0]}(:=f)$, $f^{[1]}$ и $f^{[2]}$ порождены матрицей $F$.

Важным в спектральной теории операторов Штурма-Лиувилля является вопрос об определении дефектных чисел оператора $L_0$ при данных матриц-функциях $P,\,Q$ и $R$, в частности, перечисление условий на эти матриц-функции, которые обеспечивают реализацию заданной пары $(d_+,d_-)$.
Начиная с работы В.Б. Лидского \cite{Lidskii} этот вопрос для операторов, порождённых выражением вида $-(Pf^{\prime})^{\prime}+Qf$, где $P(x)$ и $Q(x)$ - эрмитовы матрицы порядка $n$ таковы, что элементы матриц $P^{-1}(x),\,Q(x)$ локально интегрируемы на $R_+$, находится в центре внимания многих математиков (см. работы \cite{Anderson}-\cite{Ser} и списки цитированной в них литературы).

С другой стороны, класс операторов, порождённых выражением (\ref{trivial}), намного шире и, в частности, охватывает часто встречающиеся в математической физике операторы, порождённые некоторыми дифференциальными выражениями с коэффициентами-распределениями (см. ниже, п. 1.2). Основная цель настоящей работы --- изучение вопроса об индексе дефекта оператора $L_0$, порождённого выражением $l[f]$, в терминах матриц-функций $P,\,Q$ и $R$. Различные аспекты этого вопроса обсуждались также в работах \cite{Safonova8}-\cite{BMS}.

Справедлива следующая теорема.
\begin{theorem}
 Для того чтобы дефектные числа оператора $L_0$ были не максимальны, необходимо и достаточно, чтобы для некоторой последовательности не-\\пересекающихся интервалов $(a_k,b_k)\subset R_+$, $k=1,2,\ldots$, выполнялось условие
\begin{equation}
\label{80}
\sum\limits_{k=1}^{+\infty}\left (\int\limits_{a_k}^{b_k}dx
\int\limits_{a_k}^{x}||K(x,t)||^2dt\right )^{1/2}=\infty,
\end{equation}

где $K(x,t)$ --- функция Коши уравнения $l[f]=0$, т.е. решение этого уравнения по переменной $x$, удовлетворяющее начальным условиям $K(x,t)|_{x=t}=O$, $K^{[1]}(x,t)|_{x=t}=I$, а символ $||\cdot||$ - означает самосопряжённую матричную норму.
\end{theorem}
Справедливость  теоремы, аналогичной теореме 1, для случая скалярных дифференциальных операторов любого порядка и весовых пространств $\mathcal{L}_w^p$ $(1 \le  p < +\infty)$ установлена в \cite{M1}. Позже подробное доказательство этой теоремы для случая скалярных дифференциальных операторов второго порядка и пространств $\mathcal{L}_w^2$ было приведено в \cite{M3}. Там же и в \cite{M4} для случая $n=1$ приведены некоторые следствия теоремы 1.
Кроме того, достаточность условия (\ref{80}) для не максимальности дефектных чисел оператора $L_0$, порождённого выражением вида $-(Pf^{\prime})^{\prime}+Qf$ (см. выше),  используется в работе \cite{S1}. Однако, полное доказательство теоремы 1 нигде не опубликовано. Мы в данной работе устраняем этот пробел и приводим его в начале параграфа 2.

{1.2.}
Пусть $P_0, Q_0$ и $P_1$ --- эрмитовы матриц-функции порядка $n$ с измеримыми по Лебегу элементами такие, что $P^{-1}_0$ - существует и $||P^{-1}_0||, ||P^{-1}_0||||P_1||^2,$ $||P^{-1}_0||||Q_0||^2$ локально интегрируемы по Лебегу. Пусть  $\varphi:=P_1+iQ_0$ и $\varphi^*:=P_1-iQ_0$. Рассмотрим матрицу
$$
F=\begin{pmatrix}
P^{-1}_0\varphi&P^{-1}_0\\
-{\varphi^*}P^{-1}_0\varphi&-{\varphi^*}P^{-1}_0\\
\end{pmatrix}.
$$
Используя свойства матричных норм, эрмитовость матриц-функций $P_0,\,Q_0$ и $P_1$ и неравенство Коши-Буняковского, легко заключить, что для любого $[\alpha,\beta]\subset R_+$ справедливы неравенства
$$
\int\limits_{\alpha}^{\beta}||\varphi^* P^{-1}_0||=
\int\limits_{\alpha}^{\beta}||P^{-1}_0\varphi||
\int\limits_{\alpha}^{\beta}{\sqrt{||P^{-1}_0||}}\cdot {\sqrt{||P^{-1}_0||}}||\varphi||$$
$$\leq
\left(\int\limits_{\alpha}^{\beta}{||{P^{-1}_0}||}\right)^{1/2}\cdot
\left(\int\limits_{\alpha}^{\beta}{||{P^{-1}_0}||\cdot||\varphi||^2}\right)^{1/2}<+\infty
$$
и
$$
\int\limits_{\alpha}^{\beta}||\varphi^*P^{-1}_0\varphi||\leq
\int\limits_{\alpha}^{\beta}{||{P^{-1}_0}||\cdot||\varphi||^2}<+\infty.
$$

Таким образом, элементы всех блоков матрицы $F$ принадлежат пространству $L^1_{loc}(R_+)$.

Посредством матрицы $F$ определим квазипроизводные $f^{[0]},\,f^{[1]},\,f^{[2]}$, полагая, как и ранее,
$$f^{[0]}=f,\quad f^{[1]}=P_0f^{\prime}-\varphi f,\quad
 f^{[2]}=(f^{[1]})^{\prime}+{{\varphi^*}}{P^{-1}_0}f^{[1]}+{{\varphi^*}}{P^{-1}_0}\varphi f.$$

Далее, применяя замечание, сделанное в п. 1.1., заключаем, что для уравнения
\begin{equation*}
\label{t1}
-f^{[2]}=\lambda f
\end{equation*}
 справедлива теорема существования и единственности решения задачи Коши, поставленная в произвольной точке $R_+$. Кроме того, можно показать, что условия, перечисленные выше на функции $||P^{-1}_0||,\,||Q_0||$ и $||P_1||$, являются самыми общими условиями, обеспечивающими это (см. \cite{Zettl}, Ch. 1, Th. 1.2.3).

Если предположить, что элементы матрицы $P_0$ также принадлежат  $\mathcal{L}^1_{loc}(R_+)$, то легко заметить, что и элементы матрицы $\varphi$ будут локально интегрируемы на $R_+$.
Из этих предположений можно заключить, что если $^\prime$ трактовать как операцию взятия производной в смысле теории распределений, то в  выражении $f^{[2]}$ можно раскрыть  все скобки, и для него получим  формулу
\begin{equation*}
\label{p00}
f^{[2]}=(P_0f^{\prime})^{\prime}-i((Q_0f)^{\prime}+Q_0f^{\prime})-P^{\prime}_1f.
\end{equation*}
Особо отметим, что в этом равенстве координаты вектор-функций $(P_0f^{\prime})^{\prime}$, $(Q_0f)^{\prime}$ и $P^{\prime}_1f$ являются сингулярными обобщёнными функциями, причём координаты первых двух из них ---  производные в смысле теории распределений регулярных обобщенных функций, а координаты вектор-функций $Q_0f^{\prime}$ и $f^{[2]}$ --- регулярные обобщённые функции.
Таким образом, выражение $l[f]$ (см. (\ref{trivial})) в терминах обобщённых функций формально записывается в виде
\begin{equation}
\label{222}
l[f]=-(P_0f^{\prime})^{\prime}+i((Q_0f)^{\prime}+Q_0f^{\prime})+P^{\prime}_1f.
\end{equation}

В частности, если $P(x)=I$, $R(x)=\sigma(x)$ и $Q(x)=-\sigma^2(x)$, где $\sigma(x)$ - вещественнозначная, симметрическая матриц-функция порядка $n$ такова, что элементы матрицы $\sigma^2(x)$ локально интегрируемы на $R_+$, то векторное квазидифференциальное выражение (\ref{222}) примет вид
\begin{equation}
\label{p1}
l[f]=-f^{\prime\prime}+\sigma^{\prime}f.
\end{equation}

Подробный анализ билинейной формы из формулы (\ref{1.1}) позволяет установить справедливость следующей теоремы (см., напр., \cite{Safonova8} и \cite{Safonova10}).
\begin{theorem}
Пусть существует последовательность попарно непересекающихся интервалов $(a_k,b_k)\subset R_+$ $(k=1,2,\ldots)$ такая, что
элементы матрицы $\sigma$ абсолютно непрерывны на $[a_k,b_k]$ и $\sigma^{\prime}(x)\geq O$ п.в. при $x\in[a_k,b_k]$, и пусть выполнено условие
\begin{equation*}
\label{23}
\sum\limits_{k=1}^{+\infty}(b_k-a_k)^2=+\infty.
\end{equation*}
Тогда для выражения $l[f]$ (оператора $L_0$) имеет место случай предельной точки.
\end{theorem}

Истоки этой теоремы восходят к Ф. Хартману, Р.С. Исмагилову, Ф. Аткинсону и М.С.П. Истхему. В векторном случае при условии, что элементы матрицы $\sigma^{\prime}(x)$ из $L^1_{loc}(R_+)$, аналогичная теорема была доказана В.П. Серебряковым (см., напр. \cite{S2}).

Пусть теперь $x_k$ $(k=1,2,\ldots)$ -- возрастающая последовательность положительных чисел, $x_0=0$ и  $\lim\limits_{k \to +\infty}x_k=+\infty$. Если при этом в выражении (\ref{p1}) положить $\sigma(x)=C_k$ при $x\in[x_{k-1},x_k)$, где $C_k$ - вещественная симметрическая числовая матрица, то это выражение принимает вид
\begin{equation}
\label{p2}
l[f]=-f^{\prime\prime}+
\sum\limits_{k=1}^{+\infty}
{\mathcal H}_k\delta(x-x_{k})f,
\end{equation}
где $\mathcal{H}_k=C_{k+1}-C_k$, а $\delta(x)$ --- $\delta$-функция Дирака.
Таким образом, теория операторов, порожденных выражениями вида (\ref{p2}), включается в теорию операторов, порождённых векторными квазидифференциальными выражениями второго порядка.

Пусть $l^2_n$ - гильбертово пространство бесконечных последовательностей $n$-компо-\\нентных вектор-столбцов со стандартным скалярным произведением.
Справедлива следующая теорема.
\begin{theorem}
Для выражения $l[f]$ (см. (\ref{p2})) имеет место случай предельного круга в том и только в том случае, когда все решения разностного векторного уравнения
\begin{equation*}
-\frac{Z_{k+1}}{r_{k+1}r_{k+2}d_{k+1}}+\frac{1}{r^2_{k+1}}
[{\mathcal H}_k+(\frac{1}{d_k}+\frac{1}{d_{k+1}})I]Z_k-\frac{Z_{k-1}}{r_kr_{k+1}d_k}=0,\,\,\,\,\,\,k=1,2,\ldots
\end{equation*}
принадлежат пространству $l^2_n$, где $d_k=x_{k}-x_{k-1}$, $r_{k+1}=\sqrt{d_k+d_{k+1}}$.
\end{theorem}

Эта теорема  принадлежит авторам (см., напр., \cite{Safonova8} и \cite{Safonova10}) и является обобщением одной теоремы, полученной М.М. Маламудом и А.С. Костенко для случая $n=1$ (см. \cite{KoM} и \cite{KoM1}).

 {1.3.}В связи с теоремой 3 нам понадобятся следующие сведения из теории операторов, порождённых разностными выражениями второго порядка в гильбертовом пространстве $l^2_n$, приведённые, например, в \cite{Berezanski}.

Пусть $A_j,\,B_j$ ($j=0,1,\ldots$) - квадратные матрицы порядка $n$, причём $B^{-1}_j$ существуют, а $A_j$ - самосопряжены. Рассмотрим векторное разностное выражение второго порядка
\begin{equation}
\label{rv}
({\l}u)_j=B_ju_{j+1}+A_ju_j+B^*_{j-1}u_{j-1},\,\,\,\,\,j=1,2,\ldots,
\end{equation}
где $u_0,u_1,\ldots\in C^n$,
и обобщённую якобиеву матрицу
\begin{equation*}
\label{J}
J=\begin{pmatrix}
A_0&B_0&O&O&\ldots\\
B^*_0&A_1&B_1&O&\ldots\\
O&B^*_1&A_2&B_2&\ldots\\
\vdots&\vdots&\vdots&\vdots&\ddots\\
\end{pmatrix}.
\end{equation*}
Выражение (\ref{rv}) и граничные условия $A_0u_0+B_0u_1=0$ определяют в пространстве $l^2_n$ минимальный замкнутый симметрический оператор ${\L}_0$ с всюду плотной областью определения и индексом дефекта $(d_+,d_-)$ ($0\leq d_+,d_-\leq n$ и если одно из них равно $n$, то и другое такое же). Согласно терминологии, принятой в матричной проблеме моментов, говорят, что оператор ${\L}_0$ порождён матрицей $J$. Кроме того, если $d_+=d_-=0$, то говорят, что для оператора ${\L}_0$ имеет место определённый, а если  $d_+=d_-=n$, то вполне неопределённый случаи (при $n=1$ слово "вполне" опускают). Определённый случай соответствует случаю предельной точки, а вполне неопределённый - случаю предельного круга для оператора Штурма-Лиувилля.

Теорема 3 утверждает, что дефектное число оператора ${L}_0$, порождённого выражением (\ref{p2}) в гильбертовом пространстве $\mathcal{L}^2_n(R_+)$, равно $2n$ в том и только том случае, когда для матрицы $J$ имеет место вполне неопределённый случай, где элементы этой матрицы определяются формулами
\begin{equation}
\label{matrix}
A_k=\frac{1}{r^2_{k+1}}[\mathcal{H}_k+(\frac{1}{d_k}+\frac{1}{d_{k+1}})I],\,\,
B_k=-\frac{1}{r_{k+1}r_{k+2}d_{k+1}}I,\,\,\,\,\,\,\,k=1,2,\ldots,
\end{equation}
а $A_0,\,B_0$ - произвольные вещественные симметрические матрицы порядка $n$ такие, что $B^{-1}_0$ существует.

Имеет место следующая теорема - аналог теоремы 1 для операторов, порождённых обобщёнными якобиевыми матрицами.
\begin{theorem}
Для того, чтобы для оператора ${\L}_0$ имел место вполне неопределённый случай, необходимо и достаточно, чтобы для любой последовательности отрезков натуральных чисел $[n_k,m_k]$ таких, что $m_k\leq n_{k+1}\leq m_{k+1}$ ($k=1,2,\ldots$), выполнялось условие
\begin{equation*}
\label{koshi}
\sum\limits^{+\infty}_{k=1}\left(\sum\limits^{m_k}_{i=n_k}\sum\limits^{i}_{j=n_k}||K_{ij}||^2\right)^{1/2}<+\infty,
\end{equation*}

где $K_{ij}$ - решение разностного уравнения $(lu)_j=0$ с начальными условиями $K_{jj}=O$ и $K_{j+1,j}=B^{-1}_j$ ($i,j=1,2\ldots$).
\end{theorem}

В работах \cite{KM}-\cite{KM2} обсуждаются вопросы, связанные с индексами дефекта операторов, порождённых обобщёнными якобиевыми матрицами, в частности, в \cite{KM} установлена справедливость теоремы 4.

{\bf 2.  Теоремы о не максимальности и минимальности дефектных чисел}
{2.1.}
{ \large Доказательство теоремы 1.}
Пусть квадратные матриц-функции $\Phi$ и $\Psi$ порядка $n$ являются матричными решениями уравнения $l[f]=0$, удовлетворяющими начальным условиям $\Phi(0)=\Psi^{[1]}(0)=I$ и $\Phi^{[1]}(0)=\Psi(0)=O$. Сначала покажем, что матриц-функции $K(x,t),\,\Phi(x)$ и $\Psi(x)$ связаны соотношением
\begin{equation}
\label{77}
K(x,t)=\Psi(x)\Phi^*(t)-\Phi(x)\Psi^*(t).
\end{equation}

Действительно, из определения матриц-функций $\Phi$ и $\Psi$ легко извлечь, что вектор-\\столбцы матрицы
\begin{equation*}
\label{2}
T=\begin{pmatrix}
\Phi & \Psi\\
\Phi^{[1]}&\Psi^{[1]}\\
\end{pmatrix}
\end{equation*}
образуют фундаментальную систему решений системы (\ref{ss}) при $\lambda=0$. Кроме того, легко уствновить, что
$$T^{-1}=\begin{pmatrix}
\Psi^{*[1]} & -\Psi^*\\
-\Phi^{*[1]}&\Phi^*\\
\end{pmatrix}$$
(см. замечание 1 из \cite{Safonova10}).

Пусть матриц-функция $F$ определена в формуле (\ref{ss}) и вектор-функция $H:=(0,h)^t$, где $h$ - $n$-компонентный вектор-столбец с локально интегрируемыми координатами. Применяя метод вариации постоянных к неоднородной системе дифференциальных уравнений первого порядка $Y^{\prime}=FY+H$, получим, что общее решение этой системы имеет вид
$$Y=T(x)
\begin{pmatrix}
C_1  \\
C_2\\
\end{pmatrix}+
\int\limits_{x_0}^{x}
T(x)T^{-1}(t)
\begin{pmatrix}
0  \\
h \\
\end{pmatrix}dt,$$
где $C_1$ и $C_2$ - произвольные постоянные $n$-компонентные вектор-столбцы, т.е. функция Коши этой системы равна $T(x)\cdot T^{-1}(t)$.
 Из последнего равенства заключаем, что для общего решения уравнения $l[f]= h$ справедлива равенство
 $$f=\Phi C_1+\Psi C_2+
 \int\limits_{x_0}^{x}
\{\Psi(x)\Phi^*(t)-\Phi(x)\Psi^*(t)\}h(t)
dt,$$
т.е. для $K(x,t)$ справедливо равенство (\ref{77}). При этом справедливость соотношений $K(x,t)|_{x=t}=O$ и $K^{[1]}(x,t)|_{x=t}=I$ немедленно следует из представления матриц-функции $T$ и равенства $T(x)T^{-1}(x)=I$ (здесь $I$ - единичная матрица порядка $2n$).

Д о с т а т о ч н о с т ь. Из формулы (\ref{77}) и из самосопряжённости матричной нормы следует, что
$$||K(x,t)||\leq ||\Psi(x)||||\Phi(t)||+||\Phi(x)||||\Psi(t)||.$$
Записав правую часть этого неравенства в виде скалярного произведения двумерных векторов $(||\Psi(x)||,||\Phi(x)||)$ и $(||\Phi(t)||,||\Psi(t)||)$ и применяя неравенство Коши-Буня-\\ковского, находим, что
$$||K(x,t)||^2\leq (||\Phi(x)||^2+||\Psi(x)||^2)(||\Phi(t)||^2+||\Psi(t)||^2).$$

Интегрируя последнее неравенство по $t$ в пределах от $a$ до $x$, а затем по $x$ в пределах от $a$ до $b$, где $0\leq a\leq b<+\infty$, после элементарных вычислений находим, что
\begin{equation}
\label{78}
\int\limits_{a}^{b}(||\Phi(x)||^2+||\Psi(x)||^2)dx\geq \sqrt{2}\left\{\int\limits_{a}^{b}dx
\int\limits_{a}^{x}||K(x,t)||^2dt\right\}^{1/2}.
\end{equation}

Пусть теперь $(a_k,b_k)\subset R_+$ - последовательность непересекающихся интервалов, удовлетворяющая условию (\ref{80}). Применяя неравенство
\begin{equation*}
\label{81}
\int\limits_{0}^{+\infty}(||\Phi||^2+||\Psi||^2)\geq \sum\limits_{k=1}^{\infty}\int\limits_{a_k}^{b_k}(||\Phi||^2+||\Psi||^2),
\end{equation*}
а затем неравенство (\ref{78}) к каждому слагаемому из суммы, стоящей в правой части, и равенство (\ref{80}), находим, что
\begin{equation}
\label{555}
\int\limits_{0}^{+\infty}(||\Phi||^2+||\Psi||^2) =+\infty.
\end{equation}
Таким образом, либо матриц-функция $\Phi$, либо $\Psi$ содержит вектор-столбец,
не принадлежащий пространству $\mathcal{L}^2_n(R_+)$, т.е. при некотором $\lambda$ (а именно при $\lambda=0$) не все решения уравнения (\ref{000}) принадлежат этому пространству. Следовательно, дефектные числа оператора  $L_0$ не максимальны (см. п. 1.1). Достаточность доказана.

Н е о б х о д и м о с т ь. Предположим, что для оператора $L_0$ не реализуется случай предельного круга, т.е. выполняется равенство (\ref{555}). Покажем, что тогда для любого $a\in R_+$ выполняется равенство
\begin{equation}
\label{777}
\int\limits_{a}^{+\infty}dx
\int\limits_{a}^{x}||K(x,t)||^2dt=+\infty.
\end{equation}
Допустим, что это не так, т.е. для некоторого $a_0\in R_+$
\begin{equation*}
J:=\int\limits_{a_0}^{+\infty}dx
\int\limits_{a_0}^{x}||K(x,t)||^2dt<+\infty.
\end{equation*}
Применяя равенство (\ref{77}), заметим, что тогда
\begin{equation*}
\int\limits_{a_0}^{+\infty}dx
\int\limits_{a_0}^{+\infty}||K(x,t)||^2dt=2J<+\infty.
\end{equation*}
Следовательно, согласно теореме Фубини, каждый столбец матриц-функции $K(\cdot,t)$ для п.в. $t\geq a_0$ принадлежит пространству $\mathcal{L}^2_n(R_+)$. Теперь заметим, что, во-первых, существуют числа $t^0_1$ и $t^0_2$ ($> a_0$), такие что
\begin{equation*}
\det \begin{pmatrix}
\Phi(t_1) & \Psi(t_1)\\
\Phi(t_2) & \Psi(t_2)\\
\end{pmatrix} \neq 0
\end{equation*}
при $t_1=t^0_1$ и $t_2=t^0_2$,
поскольку противоположное допущение
приводит к противоречивому равенству
$\det T(t^0_1)=0$
при некотором $t^0_1> a_0$, а, во-вторых, в силу непрерывности матриц-функций $\Phi(t)$ и $\Psi(t)$ при некотором $\delta>0$ это неравенство справедливо при всех $t_1\in(t^0_1-\delta,t^0_1+\delta)$ и $t_2\in(t^0_2-\delta,t^0_2+\delta)$. Таким образом, существуют числа $t^1_1$ и $t^1_2$ из соответствующих интервалов такие, что каждый столбец матриц-функции $K(\cdot,t^1_i)$ принадлежит пространству $\mathcal{L}^2_n(R_+)$ при $i=1,2$. Кроме того, вектор-столбцы матриц-функции $K(\cdot,t^1_1)$ и $K(\cdot,t^1_2)$ линейно независимы. Действительно, из равенства
$$K(x,t^1_1)\alpha+K(x,t^1_2)\beta=0,$$
где $\alpha=(\alpha_1,\ldots,\alpha_n)^t$, $\beta=(\beta_1,\ldots,\beta_n)^t$, формулы (\ref{77}) и начальных условий $\Phi(0)=\Psi^{[1]}(0)=I$ и $\Phi^{[1]}(0)=\Psi(0)=O$ следует, что вектора $\alpha$ и $\beta$ удовлетворяют однородной системе линейных уравнений
$$\begin{cases}
\Phi^*(t^1_1)\alpha+\Phi^*(t^1_2)\beta=0\\
\Psi^*(t^1_1)\alpha+\Psi^*(t^1_2)\beta=0\\
\end{cases},$$
определитель которой отличен от 0. Следовательно, $\alpha_1=\ldots=\alpha_n=\beta_1=\ldots=\beta_n=0$.

Таким образом, все решения уравнения (\ref{000}) при $\lambda=0$ принадлежат пространству $\mathcal{L}^2_n(R_+)$, а это противоречит условию (\ref{555}). Равенство (\ref{777}) доказано. Остаётся заметить, что если выполняется это равенство, то можно найти последовательность непересекающихся интервалов $(a_k,b_k)\subset R_+,\,k=1,2,\ldots$, например, таких, что
\begin{equation*}
\int\limits_{a_k}^{b_k}dx
\int\limits_{a_k}^{x}||K(x,t)||^2dt\geq 1,
\end{equation*}
а для этой последовательности интервалов равенство (\ref{80}), очевидно, выполняется. Теорема 1 доказана.

Пусть $(a_k,b_k)\subset R_+$ ($k=1,2,\ldots$) - последовательность непересекающихся интервалов. Заметим что, как известно, при фиксированном $k$ в треугольнике $\{(x,t)| \, a_k<t \le x< b_k \}$ функция $K(x,t)$ однозначно определяется значениями коэффициентов выражения $l$. Таким образом, если выполняется условие (\ref{80}), то независимо от значений этих функций вне множества $U_{k=1}^{+\infty}{[a_k,b_k]}$ для выражения $l$ не реализуется случай предельного круга, требуется лишь, чтобы матриц-функции $P(x), Q(x)$ и $R(x)$ на $R_+$ удовлетворяли условиям, перечисленным в п.1.1. В дальнейшем мы воспользуемся этим замечанием.

{2.2.}
Справедлива следующая лемма.
\begin{lemma}
Пусть $a,b$ и $c$ - положительные числа, такие что $a<c<b$, $H$ - произвольная вещественная симметрическая матрица порядка $n$ и пусть $K(x,t)=(k_{ij}(x,t))^{n}_{i,j=1}$ - функция Коши векторного квазидифференциального уравнения
\begin{equation}
\label{1}
-y^{\prime\prime}+H\delta(x-c)y=0.
\end{equation}
Тогда
$$\int\limits_{a}^{b}dx\int\limits_{a}^{x}|k_{ii}(x,t)|^2dt
\geq\frac{1}{3\sqrt{3}}(\rho s)^2(\rho+s)\left|h_{ii}+\frac{3}{2}\left (\frac{1}{\rho}+\frac{1}{s}\right )\right|
$$
и
$$\int\limits_{a}^{b}dx\int\limits_{a}^{x}|k_{ij}(x,t)|^2dt
=\frac{|h_{ij}|^2}{9}(\rho s)^3,\mbox{если}\,\,\,\,\,\,\,i\neq j,
$$
где $\rho=c-a$ и $s=b-c$.
\end{lemma}
{\large Доказательство}. В определении квазидифференциального выражения $l[f]$ по формуле (\ref{p1}) положим
$$
\sigma(x)=\left\{
\begin{array}{rcl}
O,\mbox{если}\,\,\,a\leq x<c\\
H,\mbox{если}\,\,\,c\leq x \leq b\\
\end{array}
\right..
$$
 Тогда уравнение $l[f]=0$ примет вид (\ref{1}).
Пусть $t\in[a,b]$ фиксировано. Функция $K(x,t)$ удовлетворяет уравнению (\ref{1}) при $x>t$ и начальным условиям $K(x,t)|_{x=t}=O$ и  $K^{[1]}(x,t)|_{x=t}=I$. Из этого легко вывести, что
$$
K(x,t)=(x-t)I$$
при $t\in[a,c]$, $x\in[t,c]$ или при $
t\in[c,b]$, $x\in[t,b]$.

Пусть теперь
$t\in(a,c)$ и $x\in(c,b)$. Тогда $K(x,t)=C_1(t)+xC_2(t)$ при $x\in(c,b]$, где $C_1(t)$ и $C_2(t)$ - матриц-функции порядка $n$, и, следовательно,
$$K(c+0,t)=C_1(t)+cC_2(t)\,\,\mbox{и}\,\,K^{[1]}(c+0,t)=
C_2(t)-H(C_1(t)+cC_2(t)).$$
С другой стороны,
$$K(c-0,t)=(c-t)I\,\,\mbox{и}\,\,K^{[1]}(c-0,t)=I.$$
Используя далее условия непрерывности функций $K(x,t)$ и $K^{[1]}(x,t)$ в точке $x=c$, получаем, что $C_1(t)$ и $C_2(t)$ удовлетворяют системе
$$
\left\{
\begin{array}{rcl}
C_1(t)+cC_2(t)=(c-t)I\\
C_2(t)-H(C_1(t)+cC_2(t))=I
\end{array}
\right..
$$
Таким образом, $C_1(t)=-tI-(c-t)cH$ и $C_2(t)=I+(c-t)H$, и, следовательно,
$$
K(x,t)=
(x-t)I+(c-t)(x-c)H
$$
при $t\in[a,c]$, $x\in[c,b]$.

Из полученных формул для $K(x,t)$ следует, что при $1\leq i\leq n$
$$J_{ii}:=\int\limits_{a}^{b}dx\int\limits_{a}^{x}|k_{ii}(x,t)|^2dt=:J_1+2h_{ii}J_2+h_{ii}^2J_3,$$
где вычисления показывают, что
$$J_1=\frac{1}{12}(b-a)^4,\,J_2=\frac{1}{6}(b-c)^2(c-a)^2(b-a),\, J_3=\frac{1}{9}(b-c)^3(c-a)^3.$$
Таким образом,
$$J_{ii}=\frac{h_{ii}^2}{9}(\rho s)^3+\frac{h_{ii}}{3}(\rho s)^2(\rho+s)+\frac{1}{12}(\rho+s)^4.$$
Записав интеграл $J_{ii}$ в виде
$$J_{ii}=(\rho
s)^3\left[\left[\frac{h_{ii}}{3}+\frac12\left(\frac{1}{\rho}+\frac{1}{s}\right)\right]^2+\left(\frac{1}{\rho}+\frac{1}{s}\right)^2
\left[\frac{1}{12}\rho s\left(\frac{1}{\rho}+\frac{1}{s}\right)^2-\frac14\right]\right],$$
заметим, что
$$\frac{1}{12}\rho s\left(\frac{1}{\rho}+\frac{1}{s}\right)^2-\frac14\geq\frac{1}{12}.$$
Поэтому
$$J_{ii}\geq(\rho s)^3\left[\left[\frac{h_{ii}}{3}+\frac12\left(\frac{1}{\rho}+\frac{1}{s}\right)\right]^2+\frac{1}{12}\left(\frac{1}{\rho}+\frac{1}{s}\right)^2\right].$$
Применив ещё раз неравенство между арифметическим и геометрическим средними, окончательно находим, что справедливо первое утверждение леммы 2.

Кроме того, как мы уже показали выше, элемент $k_{ij}(x,t)$ матриц-функции $K(x,t)$ при $1\leq i\neq j\leq n$ определяется равенством $k_{ij}(x,t)=(c-t)(x-c)h_{ij}$. Поэтому
$$\int\limits_{a}^{b}dx\int\limits_{a}^{x}|k_{ij}(x,t)|^2dt=\frac{|h_{ij}|^2}{9}(\rho s)^3.$$
Лемма 2 доказана.

Используя лемму 2, докажем справедливость следующей теоремы.
\begin{theorem}
Пусть
$(a_k,b_k)\subset R_+$ ($k=1,2,\ldots$) - последовательность попарно непересекающихся интервалов и $c_k$ - последовательность положительных чисел, такая что $a_k<c_k<b_k$, и
на каждом отрезке $[a_k,b_k]$ выполнены равенства
$P(x)=I$, $R(x)=\sigma(x)$,  $Q(x)=-\sigma^2(x)$, где $\sigma(x)$ - кусочно-постоянная матриц-функция порядка $n$ с матрицей скачков ${\mathcal H}_k=(h^k_{ij})_{i,j=1}^{n}$ в точке $c_k$. Пусть далее, числа $\rho_k=c_k-a_k$, $s_k=b_k-c_k$ и матрица ${\mathcal H}_k$ таковы, что выполняется одно из следующих условий
\begin{equation}
\label{31}
\sum\limits_{k=1}^{+\infty}\rho_ks_k\sqrt{\rho_k+s_k}\sqrt{\left|h^k_{ii}+\frac{3}{2}\left (\frac{1}{\rho_k}+\frac{1}{s_k}\right)\right|}=+\infty
\end{equation}
хотя бы для одного $i$, такого что $1\leq i\leq n$, или
\begin{equation}
\label{32}
\sum\limits_{k=1}^{+\infty}(\rho_k s_k)^{3/2}{|h^k_{ij}|}=+\infty
\end{equation}
хотя бы для одной пары $(i,j)$, такой что $1\leq i\neq j\leq n$.\\
Тогда для оператора $L_0$ не реализуется случай предельного круга.
\end{theorem}
{\large Доказательство.}
Применим лемму 2, положив  $a=a_k$, $c=c_k$,   $b=b_k $. Тогда
$$\int\limits_{a_k}^{b_k}dx\int\limits_{a_k}^{x}|k_{ii}(x,t)|^2dt
\geq\frac{1}{3\sqrt{3}}(\rho_k s_k)^2(\rho_k+s_k)\left|h^k_{ii}+\frac{3}{2}\left (\frac{1}{\rho_k}+\frac{1}{s_k}\right )\right|
$$
и
$$\int\limits_{a_k}^{b_k}dx\int\limits_{a_k}^{x}|k_{ij}(x,t)|^2dt
=\frac{|h^k_{ij}|^2}{9}(\rho_k s_k)^3,\mbox{если}\,\,\,i\neq j,\,\,\,k=1,2,\ldots
$$

Полагая $||K(x,t)||^2=\sum\limits_{i,j=1}^{n}|k_{ij}(x,t)|^2$, и при фиксированном $1\leq i\leq n $ используя неравенство $||K(x,t)||\geq|k_{ii}(x,t)|$, получим
\begin{equation*}
\int\limits_{a_k}^{b_k}dx
\int\limits_{a_k}^{x}||K(x,t)||^2dt\geq
3^{-3/2}(\rho_k s_k)^2(\rho_k+s_k)\left|h^k_{ii}+\frac{3}{2}\left (\frac{1}{\rho_k}+\frac{1}{s_k}\right )\right|.
\end{equation*}
Извлекая квадратный корень из обеих частей этого неравенства, а затем суммируя по $k$, находим, что
\begin{equation*}
\label{00}
\sum\limits_{k=1}^{+\infty}\left (\int\limits_{a_k}^{b_k}dx
\int\limits_{a_k}^{x}||K(x,t)||^2dt\right )^{1/2}  \geq {3^{-3/4}}  \sum\limits_{k=1}^{+\infty} \rho_k s_k\sqrt{(\rho_k+s_k)}\sqrt{\left|h^k_{ii}+ \frac{3}{2}\left (\frac{1}{\rho_k}+\frac{1}{s_k}\right)\right|}.
\end{equation*}
Рассуждая аналогично с учётом второго утверждения леммы 2, для $1\leq i\neq j\leq n$ находим, что
\begin{equation*}
\label{01}
\sum\limits_{k=1}^{+\infty}\left (\int\limits_{a_k}^{b_k}dx
\int\limits_{a_k}^{x}||K(x,t)||^2dt\right )^{1/2}  \geq
\frac 13\sum\limits_{k=1}^{+\infty}(\rho_k s_k)^{3/2}{|h^k_{ij}|}.
\end{equation*}
Остается применить теорему 1 с учётом условий (\ref{31}) и (\ref{32}). Теорема 5 доказана.

{2.3.}
 Приведем некоторые следствия из теоремы 5.

Формулировка теоремы 5, очевидно, упрощается, если предположить, что точки $c_k$ являются серединами отрезков $[a_k,b_k]$ ($k=1,2,\ldots$). В этой ситуации справедливо

\begin{corollary}
Пусть выполнены условия теоремы 5, $c_k=\frac{a_k+b_k}{2}$ и
$$\sum\limits_{k=1}^{+\infty}\rho^{5/2}_k\sqrt{\left|  {h^k_{ii}} +\frac{6}{\rho_k}\right|}=+\infty$$
хотя бы для одного $i$, такого что $1\leq i\leq n$, или
$$\sum\limits_{k=1}^{+\infty}\rho^3_k{|h^k_{ij}|}=+\infty$$
хотя бы для одной пары $(i,j)$, такой что $1\leq i\neq j\leq n$, где $\rho_k=b_k-a_k$.\\
Тогда  для оператора $L_0$ не реализуется случай предельного круга.
\end{corollary}

Применим теорему 5 в ситуации, когда дифференциальное выражение $l$ имеет вид (\ref{p2}). В этом случае  в качестве точек $c_k$ ($k=1,2,\ldots$) выберем точки $x_k$, а в качестве отрезков $[a_k,b_k]$ -- отрезки $[x_k-\frac{d_k}{2},x_k+\frac{d_{k+1}}{2}]$, где, напомним, что $d_k=x_{k}-x_{k-1}$ (см. формулировку теоремы 3). Справедливо следующее следствие теоремы 5.

\begin{corollary}
Пусть дифференциальное выражение $l[f]$ имеет вид (\ref{p2}) и
$$\sum\limits_{k=1}^{+\infty}d_kd_{k+1}r_{k+1}\sqrt{\left| {h^k_{ii}} +\frac{3}{2}\left (\frac{1}{d_k}+\frac{1}{d_{k+1}}\right) \right|}=+\infty$$
хотя бы для одного $i$, такого что $1\leq i\leq n$, или
$$\sum\limits_{k=1}^{+\infty}(d_kd_{k+1})^{3/2}|h^k_{ij}|=+\infty$$
хотя бы для одной пары $(i,j)$, такой что $1\leq i\neq j\leq n$.\\
Тогда для оператора $L_0$ не реализуется случай предельного круга.
\end{corollary}

{2.4.}
В пп. 2.2, 2.3 настоящего параграфа было продемонстрировано, как из теоремы 1 можно получить признаки не максимальности дефектных чисел для векторных операторов Штурма-Лиувилля, порожденных выражениями вида (\ref{trivial}),  (\ref{p1}) или  (\ref{p2}). Применяя теорему 3, получаем, что из справедливости следствия 2 теоремы 5 следует, что для обобщённой якобиевой матрицы $J$ с элементами вида  (\ref{matrix})  не реализуется вполне неопределённый случай.

В этом пункте получим признак реализации определенного случая для якобиевых матриц
$J$ с элементами вида  (\ref{matrix}). Здесь же отметим, что, согласно теореме 2,  он будет одновременно и признаком реализации случая предельной точки для векторного оператора Штурма-Лиувилля, порожденного выражением вида (\ref{p2}).

Справедлива следующая теорема.
\begin{theorem}
Пусть последовательность $d_k$ такова, что
\begin{equation}
\label{5}
\sum\limits_{k=1}^{+\infty}d_k^{2}=+\infty.
\end{equation}
Тогда для  якобиевой матрицы  $J$ с элементами, определяемыми равенствами (\ref{matrix}), имеет место определенный случай.
\end{theorem}
{\large Доказательство.}
Из равенства (\ref{matrix}) следует, что
$||B_k||^{-1}=\frac{1}{\sqrt{n}}{r_{k+1}r_{k+2}d_{k+1}}$. Кроме того,
$$
d^2_{k+1}\leq \sqrt{d_{k}+d_{k+1}}\sqrt{d_{k+1}+d_{k+2}}d_{k+1}=r_{k+1}r_{k+2}d_{k+1}$$
$$\leq\frac 1 2(r^2_{k+1}+r^2_{k+2})d_{k+1}=\frac 1 2(d_kd_{k+1}+2d^2_{k+1}+d_{k+1}d_{k+2})$$
$$\leq \frac{1}{4}(d^2_k+6d^2_{k+1}+d^2_{k+2}).$$
Следовательно,
$$\frac{1}{\sqrt{n}}d^2_{k+1}\leq ||B_k||^{-1}\leq \frac{1}{4\sqrt{n}}(d^2_k+6d^2_{k+1}+d^2_{k+2}).$$
 Из этого неравенства и условия (\ref{5}) теоремы 6 следует, что
 $$\sum\limits_{k=0}^{+\infty}||B_k||^{-1}=+\infty.$$
 Последнее соотношение является условием Карлемана, обеспечивающим реализацию определенного случая для матрицы  $J$ (см., напр., \cite{Berezanski}, гл.VII, $\S$2, п.11, теорема 2.9).
 Теорема 6 доказана.

  Отметим, что, рассуждая также, как в \cite{KM}, из теоремы 4 можно получить достаточные условия реализации вполне неопределённого и не вполне неопределённого случаев для матрицы $J$. При $n=1$ эти признаки будут близки к некоторым результатам работ \cite{KoM} и \cite{KoM1}. Однако мы здесь этого делать не будем.

{\bf 3. Заключительные замечания. Примеры}

В теоремах 2 и 5 данной работы предполагается, что квазидифференциальное выражение $l$ имеет вид (\ref{p1}) на последовательности непересекающихся интервалов $(a_k,b_k)$ ($k=1,2,\ldots$), где функция $\sigma(x)$ абсолютно непрерывна на каждом отрезке $[a_k,b_k]$ в случае теоремы 2, и является ступенчатой функцией с одним скачком в случае теоремы 5, и удовлетворяет некоторым дополнительным условиям при $x\in[a_k,b_k]$, обеспечивающим реализацию случая предельной точки для выражений вида (\ref{trivial}) в случае теоремы 2 и не максимальность дефектных чисел - в случае теоремы 5. Таким образом, в этих теоремах допускается, что  коэффициенты выражения $l$ могут быть произвольными матриц-функциями вне  отрезков $[a_k,b_k]$, требуется лишь выполнение условий, перечисленных в начале п.1.1, т.е. элементы этих матриц-функций должны быть измеримы на $R_+$ и суммируемы на каждом её замкнутом конечном интервале. В частности, выражение $l$ может иметь вид (\ref{222}) при соответствующих ограничениях на коэффициенты $P_0$, $Q_0$ и $P_1$.

В работе \cite{SS} Ш. Крист и Г. Штольц показали, что если $d_j=\frac 1j$ и $\mathcal{H}_j=(-2j-1)I$ ($j=1,2,\ldots$), то для выражения $l$ вида (\ref{p2}) при $n=1$ реализуется случай предельного круга и, по-видимому, впервые обнаружили, что для этого выражения такое возможно. Позже в работах \cite{KoM}, \cite{KoM1} А.С. Костенко и М.М. Маламуда и в работе \cite{Konechnaya} Н.Н. Конечной были построены многочисленные примеры случаев реализации предельной точки или предельного круга для выражений вида (\ref{p2}) при $n=1$.

Обобщённые якобиевы матрицы вида $J$ возникают в связи с матричной степенной проблемой моментов (см., напр.,  \cite{Krein1}) и хорошо изучены. В частности, в работах \cite{KM} - \cite{KM2} установлены критерии максимальности дефектных чисел и различные признаки реализации случаев максимальности и не максимальности дефектных чисел в терминах элементов матрицы $J$. Применив эти признаки и теорему 3, можно получить условия максимальности и не максимальности дефектных чисел оператора $L_0$, порождённого выражением (\ref{p2}), в терминах ${\mathcal H}_k$ и $d_k$. В частности, используя теорему 1 из \cite{KM2} и теорему 3 данной работы, можно показать, что справедливы следующая теорема и следствие из неё.
\begin{theorem}
Пусть элементы матрицы $J$ таковы, что\\
$(a)$
$\sum\limits_{j=1}^{+\infty}\left(r_{2j+s}\frac{d_{1+s}d_{3+s}\cdot\ldots\cdot d_{2j-1+s}}{d_sd_{2+s}\cdot\ldots\cdot d_{2j-2+s}}\right)^2<+\infty,\,\,\,\,\,s=1,2,$\\
$(b)$
$\sum\limits_{j=1}^{+\infty}\left(\frac{d_{1+s}d_{3+s}\cdot\ldots\cdot d_{2j-1+s}}{d_sd_{2+s}\cdot\ldots\cdot d_{2j-2+s}}\right)^2||\mathcal{H}_{2j+s-1}+
\left(\frac{1}{d_{2j+s-1}}+\frac{1}{d_{2j+s}}\right)I||<+\infty,\,\,\,\,\,s=1,2.$\\
\\
Тогда для оператора $L_0$ реализуется случай предельного круга.
\end{theorem}
\begin{corollary}
Пусть элементы матрицы $J$ таковы, что\\
$1)$ ${r_kr_{k+3}d_kd_{k+2}}\geq {r_{k+1}r_{k+2}d^2_{k+1}}$ или
${r_kr_{k+3}d_kd_{k+2}}\leq {r_{k+1}r_{k+2}d^2_{k+1}}$ для всех $k=1,2,...,$\\
$2)$ $\sum\limits_{k=1}^{+\infty}d^2_k<+\infty$,\\
$3)$ $\sum\limits_{k=1}^{+\infty}d_{k+1}||\mathcal{H}_k+\left(\frac{1}{d_k}+\frac{1}{d_{k+1}}\right)I||<+\infty$.\\
Тогда для выражения оператора $L_0$ реализуется случай предельного круга.
\end{corollary}
Теорема 7 и следствие 3 являются обобщениями некоторых результатов из работ \cite{Safonova8} и \cite{Safonova10}.


\noindent Мирзоев К.А., \\
МГУ имени М.В. Ломоносова, \\
Ленинские Горы, 1,\\ 
119991, г. Москва, Россия, \\
email:{mirzoev.karahan@mail.ru}\\
Сафонова Т.А., \\
САФУ имени М.В. Ломоносова, \\
Набережная Северной Двины, 17,\\ 
163002, г. Архангельск, Россия, \\
email:{tanya.strelkova@rambler.ru} \\
\end{document}